\documentclass[a4paper,12pt,fleqn]{article}
\usepackage{times}
\usepackage{amsthm, amssymb, amsmath}
\usepackage{epsfig}

\newcommand{\pl}{Pochhammer}

\newcommand{\aph}{\alpha}

\newcommand{\bgm}[1]{\Gamma\lrdl{(}{)}{#1}}

\newcommand{\lrdl}[3]{\left #1 {#3} \right #2}

\newcommand{\qhl}{$q$-hypergeometric}

\newcommand{\id}{identities}

\newcommand{\rr}{Rogers-Ramanujan}
\newcommand{\slz}{Saalsch\"{u}tz}

\newcommand{\bib}{\bibitem }

\newcommand{\z}[2]{\genfrac{}{}{1pt}{0}{#1}{#2}}

\newcommand{\efrc}[3]{#1_{r} = \z{#2}{#3}}
\newcommand{\ez}[3]{#1_{r} = \z{#2}{#3}}

\newcommand{\efnc}[3]{#1_{n} = \z{#2}{#3}}

\newcommand{\sff}[2]{\lrdl{(}{)}{#1}_{#2}}

\newcommand{\pfq}[5]{ _{#1}F_{#2} \left[ \begin{array}{c}#3\\[9pt] #4 %%@
\end{array}{;#5} \right ]}

\newcommand{\smm}{\sum^{\infty}_{n=0}}

\newcommand{\suu}[2]{\sum^{#2}_{#1}}

\newcommand{\dixon}{ ~\cite[ p. 243, (III. 8)]{slt3} }
\newcommand{\nfff}{ ~\cite[ p. 244, (III. 12)]{slt3} }
\newcommand{\fff}{ ~\cite[ p. 244, (III. 13)]{slt3} }
\newcommand{\nfft}{ ~\cite[ p. 245, (III. 22)]{slt3} }

\newcommand{\nmfft}{ ~\cite[ p. 244, (III. 11)]{slt3} }
\newcommand{\dougall}{ ~\cite[ p. 244, (III. 14)]{slt3} }

\newcommand{\lpfq}[4]{ _{#1}F_{#2} \left[ \begin{array}{c}#3\\[9pt] #4 %%@
\end{array}\right. }
\newcommand{\rpfq}[3]{\left. \begin{array}{c}#1\\[9pt] #2 \end{array};#3 %%@
\right]}
\newcommand{\beqn}[2]{\begin{equation*}\hspace{-1cm} #1 \tag{#2}
\end{equation*}} 
\newcommand{\beq}[1]{\begin{equation*} #1 \end{equation*}}
\newcommand{\bmltl}[2]{\begin{\multline} #1 \tag{#2}\end{multline}} 
\title{\bf Extensions of the classical theorems for very well-poised hypergeometric functions\footnote{``This is a pre-print of (a part of) an article published in Revista de la Real Academia de Ciencias Exactas, F\'{i}sicas y Naturales. Serie A. Matem\'{a}ticas (RACSAM). The final authenticated version is available online at: https://doi.org/10.1007/s13398-017-0485-5 ". } }
\author{\bf Yashoverdhan Vyas and Kalpana Fatawat}
\date{}
\begin{document}
\maketitle
\begin{abstract}
\noindent The classical summation and transformation theorems for very well-poised hypergeometric functions, namely, $_{5}F_4(1)$ summation, Dougall's $_{7}F_6(1)$ summation, Whipple's $_{7}F_6(1)$ to $_{4}F_3(1)$ transformation   and Bailey's $_{9}F_8(1)$ to $_{9}F_8(1)$ transformation  are extended. These extensions are derived by applying the well-known Bailey's transform method along with the classical very well-poised summation and transformation theorems for very well-poised hypergeometric functions and the Rakha and Rathie's extension of the Saalsch\"{u}tz's theorem. To show importance and applications of the discovered extensions, a number of special cases  are pointed out, which leads not only to the extensions of other classical theorems for very well-poised and well-poised hypergeometric functions but also generate new hypergeometric summations and transformations.
\end{abstract}

\section{ Introduction} 
The generalized hypergeometic function ~\cite{aa}, for any integer $r\geq 1$ and parameter vector $(a;b)=(a_{1},a_{2}, ... , a_{r+1}; b_{1}, b_{2}, ...,b_{r}) \in \mathbb{C}^{r+1}\times\mathbb{C}^{r}$ in which no $b_i$ is a nonpositive integer, is defined as:
\beq{\pfq{r+1}{r}{a_{1} , a_{2} , ... , a_{r+1}}{b_{1} , b_{2} , ... , b_{r}}{z}= \smm \ \z{\sff{a_{1}}{n} \sff{a_{2}}{n} ... \sff{a_{r+1}}{n}}{\sff{b_{1}}{n} \sff{b_{2}}{n} ... \sff{b_{r}}{n} \ n!}z^{n}}
\\ 
where $\sff{a}{n}$ stands for \pl's shifted factorial $a(a+1)(a+2)...(a+n-1)$ with $(a)_{0}=1$. The series converges absolutely on $|z| < 1$,
and if the parametric excess $s=\suu{i=1}{r}b_{i}-\suu{i=1}{r+1}a_{i}$ has positive real part, it converges at $z = 1$ and for convergence at $z = -1, Re(s) > -1$ suffices. \\ 

A generalized hypergeometric series or function is said to be $s$-balanced if the
parametric excess equals $s$ and balanced if $s=1$. It is well-poised if the parameters $ \{a_{i}\}, \{b_{i}\}$ can be separately permuted so that $1+a_{1} = a_{2} + b_{1} = a_{3} + b_{2} = ... = a_{r+1}+b_{r}$ and very well-poised if $a_{2}=1+ \frac{a_{1}}{2}$ along with the condition $1+a_{1} = a_{2} + b_{1} = a_{3} + b_{2} = ... = a_{r+1}+b_{r}$ for well-poised nature. \\ 

It is well known that whenever a generalized hypergeometric function reduces to quotient of the products of the gamma function, the results are very important from the application point of view in numerous areas of physics, mathematics and statistics including in series systems of symbolic computer algebra manipulation. There are certain classical summation and transformation theorems  ~\cite[ p. 243-245, Appendix III]{slt3}, namely, Gauss theorem for $_{2}F_1(1)$, Kummer theorem for well-poised $_{2}F_1(-1)$, Dixon theorem for well-poised $_{3}F_2(1)$, \slz \ theorem for balanced $_{3}F_2(1)$, Dougall theorem for a very well-poised $_{7}F_6(1)$, Whipple transformation  ~\cite[ p. 61, (2.4.1.1)]{slt3}  of very well-poised $_{7}F_6(1)$ to balanced $_{4}F_3(1)$, very well-poised summations for $_{5}F_4(1)$,  $_{4}F_3(\pm 1)$ and above all the Bailey's ~\cite[ p. 71, (2.4.4.1)]{slt3} very well-poised transformation of $_{9}F_8(1)$ to $_{9}F_8(1)$. The Bailey's $_{9}F_8(1)$ transformation is mentioned in Bailey ~\cite[\textsection 7.6, p.63 and \textsection 4.3(7), p.27]{bb} and Slater ~\cite[ p. 71, (2.4.4.1)]{slt3} as one of the most general known transformation of terminating very well-poised series which has been found. Also, the role of Bailey's $_{9}F_8(1)$ transformation in obtaining the two most general expansions can be found from ~\cite{y3} . Further, all the above mentioned summations and transformations can be derived from the well-known Bailey's transform, which is a very general and useful result of the theory of hypergeometric series, as shown by Slater ~\cite[ p. 58-74]{slt3}. Recently, a good progress has been achieved in generalizing various summations and transformations for hypergeometric series; for this we refer to ~\cite{kim1, kim, lavoi1, lavoi3, lavoi2, maier1, maier, miller1, miller2, miller3, miller4, prud, rakha11, 2rakha11, rathie7, vid2} and the papers cited there. \\ 

 Notably,  Maier ~\cite{maier1}, motivated from Slater's work ~\cite[ \textsection 2.6.1]{slt3}, has discovered extension of Euler's transformation  and several new hypergeometric summations with nonlinear restrictions on the parameters.  But, the extensions for the classical theorems on very well-poised $_{9}F_8(1)$, $_{7}F_6(1)$ and $_{5}F_4(1)$ hypergeoemtric functions have not appeared so far. The present research article is an attempt in this direction.\\

Following results will be needed in our work:\\

W. N. Bailey~\cite{bly1, bly2} established a series 
transform known as Bailey's transform and gave a mechanism to derive 
ordinary and \qhl\ \id\ and \rr\ type \id . Its extensions and other types of Bailey's transform have also appeared in ~\cite{yvyas, y2}.\\

It states that if  
$$ \beta_n = \sum^n_{r=0} \alpha _{r} \ u_{n-r} \ v_{n+r}\eqno (1.1)$$ and 
$$\gamma_n =\sum^\infty_{r=n} \delta_{r} \ u_{r-n} \ v_{r+n} \eqno (1.2)$$ 
 then, subject to convergence conditions,
  $$ \sum^ \infty_{n=0} \alpha_n\gamma_n = 
\sum^\infty_{n=0}\beta_n\delta_n \eqno (1.3) $$

It is outlined in Slater~\cite[p. 61]{slt3} that how the following transformation (1.4) due to Whipple, can be derived from the Bailey's transform by using \slz \ theorem and $_{5}F_4(1)$ summation:

\beq{\pfq{7}{6}{a,1+\z{a}{2}, b, c, d, e, -N}{\z{a}{2}, 1+a-b, 1+a-c, 1+a-d, 1+a-e, 1+a+N}{1}\hspace{3.5cm} }
\beqn {= \z{\sff{1+a-d-e}{N} \sff{1+a}{N}}{\sff{1+a-d}{N} \sff{1+a-e}{N}} \ \pfq{4}{3}{1+a-b-c, d, e, -N}{1+a-b, 1+a-c, d+e-a-N}{1}}{1.4}

and then by particularising the $_{4}F_{3}(1)$ to be summable by \slz \ theorem, 
 one can also obtain the Dougall's $_{7}F_6(1)$ summation.\\ 

Here, we shall also make use of the following extension of \slz \ theorem (1.5) obtained by Rakha and Rathie in ~\cite{2rakha11} :

\beqn {\pfq{4}{3}{a, b, f+1, -N}{c, f, 2+a+b-c-N}{1}= \z{\sff{c-a-1}{N} \sff{c-b-1}{N}}{\sff{c}{N} \sff{c-a-b-1}{N}} \z{\sff{g+1}{N}}{\sff{g}{N}} \ 
}{1.5} \\ 
\beq{\mbox{  where  } g=\z{f(1+b-c)(1+a-c)}{ab-f(1+a+b-c)} \hspace{11.5cm}}

Other remaining definitions and notations are from ~\cite{slt3, bb} .

\section{The extensions of classical very well-poised hypergeometric summations and transformations}

In this section, we state four new hypergeometric summations and transformations :

\beq{\pfq{7}{6}{a,1+\z{a}{2}, b, a-f+1, d, f+1, -N}{\z{a}{2}, 1+a-b, f, 1+a-d, a-f, 1+a+N}{1}\hspace{6.0cm} }
\beqn {= \z{\sff{1+a}{N}\sff{a-b-d}{N} }{\sff{1+a-b}{N} \sff{1+a-d}{N}} \z{\sff{g+1}{N}}{\sff{g}{N} }}{2.1} \\ 
\beq{\mbox{  where  } g=\z{f(b+d-a)(f-a)}{f(a-f)-db} \hspace{11.5cm}}

The summation given in (2.1) is an extension of the classical $_{5}F_4(1)$ summation theorem \fff , which follows from (2.1) by letting $f \to \infty$.

\beq{\pfq{9}{8}{f-1,\z{f+1}{2}, a_{1} , a_{2} , f-p ,p+1, d_{1}, d_{2}, -N}{\z{f-1}{2}, f-a_{1} , f-a_{2} , p ,f-p-1, f-d_{1}, f-d_{2}, f+N}{1} \hspace{3.5cm} }
\beqn {= \z{\sff{f}{N} \sff{f-d_{1}-d_{2}}{N}}{\sff{f-d_{1}}{N} \sff{f-d_{2}}{N}} \ \pfq{5}{4}{d_{1}, d_{2},f-a_{1}-a_{2}-1, h+1, -N}{f-a_{1} , f-a_{2} ,1+d_{1}+d_{2}-f-N, h}{1}}{2.2}

\beqn{\mbox{  where  } h=\z{p(a_{1}+a_{2}-f+1)(p-f+1)}{p(f-p-1)-a_{1}a_{2}}}{2.3}

The transformation (2.2) is a generalization of Whipple's result (1.4), which follows from (2.1) by letting $p \to \infty$.

\beq{\pfq{9}{8}{f-1,\z{f+1}{2}, a_{1} , f-p ,p+1, d_{1}, d_{2}, 2f-2-d_{1}-d_{2}-a_{1}+N, -N}{\z{f-1}{2}, f-a_{1} , p ,f-p-1, f-d_{1}, f-d_{2}, 2+a_{1}+d_{1}+d_{2}-f-N, f+N}{1} \hspace{3.5cm} }
\beqn {= \z{\sff{f}{N} \sff{f-d_{1}-d_{2}}{N} \sff{f-a_{1}-d_{1}-1}{N} \sff{f-a_{1}-d_{2}-1}{N}}{\sff{f-d_{1}}{N} \sff{f-d_{2}}{N} \sff{f-a_{1}}{N} \sff{f-a_{1}-d_{1}-d_{2}-1}{N}} \z{\sff{k+1}{N}}{\sff{k}{N}}}{2.4}

\beqn{\mbox{  where  } k=\z{h(1+d_{1}+a_{1}-f)(1+d_{2}+a_{1}-f)}{d_{1}d_{2}-h(1+d_{1}+d_{2}+a_{1}-f)}} {2.5}

\beqn{\mbox{  and  } h=\z{p(f-1-d_{1}-d_{2}+N)(p-f+1)}{p(f-p-1)-a_{1}(2f-2-d_{1}-d_{2}-a_{1}+N)}} {2.6}  

The summation (2.4) is a generalization of Dougall's theorem \dougall , which follows from (2.4) by selecting $p=a_1$.

\beq{\pfq{11}{10}{f-1,\z{f+1}{2}, a_{1} , f-p ,p+1, d_{1}, d_{2}, b_{1}, b_{2}, b_{3}, -N}{\z{f-1}{2}, f-a_{1} , p ,f-p-1, f-d_{1}, f-d_{2}, f-b_{1}, f-b_{2}, f-b_{3}, f+N}{1}  }
\beq {\times \z{\sff{1+c-b_{1}}{N} \sff{1+c-b_{2}}{N} \sff{1+c-b_{3}}{N} \sff{1+c-b_{1}-b_{2}-b_{3}}{N}}{\sff{1+c}{N} \sff{1+c-b_{1}-b_{2}}{N} \sff{1+c-b_{1}-b_{3}}{N} \sff{1+c-b_{2}-b_{3}}{N}}}
\beq{= \ \pfq{9}{8}{c, 1+\z{c}{2}, b_{1}, b_{2}, b_{3}, f-a_{1}-d_{2}-1, f-a_{1}-d_{1}-1, f-d_{1}-d_{2}, -N }{\z{c}{2}, 1+c-b_{1}, 1+c-b_{2}, 1+c-b_{3}, f-d_{1}, f-d_{2}, f-a_{1}, 1+c+N}{1} \hspace{2.5 cm}}
\beq{+ \z{A'B}{A D} \ \lpfq{10}{9}{c+1, 2+\z{c}{2}, b_{1}+1, b_{2}+1, b_{3}+1, f-a_{1}-d_{2}, }{1+\z{c}{2}, 2+c-b_{1}, 2+c-b_{2}, 2+c-b_{3}, 1+f-d_{1}, }}
\beq{\rpfq{f-a_{1}-d_{1}, A+1, B+2, -N+1}{1+f-d_{2}, 1+f-a_{1}, B+1, 2+c+N}{1} } 
\beq{= \ \lpfq{11}{10}{c, 1+\z{c}{2}, b_{1}, b_{2}, b_{3}, f-a_{1}-d_{2}-1, f-a_{1}-d_{1}-1, }{\z{c}{2}, 1+c-b_{1}, 1+c-b_{2}, 1+c-b_{3}, f-d_{1}, f-d_{2},} \hspace{7.0cm}}
\beqn{\rpfq{f-d_{1}-d_{2}-1, 1+\z{c}{2}+\lambda, 1+\z{c}{2}-\lambda, -N }{ f-a_{1}, \z{c}{2}-\lambda, \z{c}{2}+\lambda, 1+c+N}{1}}{2.7}

where $c=2f-2-d_{1}-d_{2}-a_{1}$, $3f = 2+b_{1}+b_{2}+b_{3}-d_{1}-d_{2}-a_{1}-N$, $A=f-d_{1}-d_{2}-1$ and $\lambda^{2} = \z{c^{2}}{4}-AD$ . The values of $B$, $D$ and $A'$   are as given below:

\beqn{B=\z{d_{1} d_{2}(p(1-f+p) + a_{1} c) + pA(1-f+p)(1+a_{1}+d_{1}+d_{2}-f)} {p(1-f+p)(1+a_{1}+d_{1}+d_{2}-f)+d_{1} d_{2}a_{1}}}{2.8}

\beqn{D=\z{-p(1-f+p)(1+a_{1}+d_{1}-f)(1+a_{1}+d_{2}-f)}{p(1-f+p) (1+a_{1}+d_{1}+d_{2}-f)+d_{1} d_{2}a_{1}}}{2.9}

\beqn{A' =\z{c (1+\z{c}{2}) b_{1} b_{2} b_{3} (f-a_{1}-d_{2}-1) (f-a_{1}-d_{1}-1)  (A) (B+1) (-N)} {(\z{c}{2}) (1+c-b_{1}) (1+c-b_{2}) (1+c-b_{3}) (f-d_{1}) (f-d_{2}) (f-a_{1})  (B) (1+c+N)}}{2.10}

The transformation (2.7) provides two forms of the generalization of Bailey's $_{9}F_8(1)$ transformation
~\cite[ p. 71, (2.4.4.1)]{slt3}. The first form of the generalization follows from the equality of first and middle terms in (2.7) and the second form of the generalization follows from the equality of first and last terms in (2.7).

The Bailey's $_{9}F_8(1)$ transformation follows from the first form of the generalization by selecting $p=a_1$ as explained below:

For $p=a_1$, $D$ becomes $f-a_{1}-1$ and $B=A$ and the last $_{10}F_{9}$ in the middle term of (2.7) becomes $_{9}F_8$. Now, using the following consequence of the last result out of the four results for $_{p}F_q$ by Rainville ~\cite[Eq.(11), p.81]{rain}:

\beq{\pfq{p}{q}{a_{1} , a_{2} , ... , a_{p}}{b_{1} , b_{2} , ..., b_{k}-1 , ... , b_{q}}{z}= \hspace{12cm}}
\beqn{\pfq{p}{q}{a_{1} , a_{2} , ... , a_{p}}{b_{1} , b_{2} , ... , b_{q}}{z}+\z{1}{b_{k}-1}\z{a_{1}  a_{2}  ...  a_{p} z}{b_{1}  b_{2}  ...  b_{q}} \  \pfq{p}{q}{a_{1}+1 , a_{2}+1 , ... , a_{p}+1}{b_{1}+1 , b_{2}+1 , ... , b_{q}+1}{z}}{2.11}

\noindent with $z=1$, $p=9$, $q=8$ and $b_{k}=f-a_{1}$, the sum of two $_{9}F_8$'s in the middle term of (2.7) converts into one $_{9}F_8$ to give

\beq{\pfq{9}{8}{f-1,\z{f+1}{2}, a_{1}+1 , d_{1}, d_{2}, b_{1}, b_{2}, b_{3}, -N}{\z{f-1}{2}, f-a_{1}-1,  f-d_{1}, f-d_{2}, f-b_{1}, f-b_{2}, f-b_{3}, f+N}{1}  }
\beq {\times \z{\sff{1+c-b_{1}}{N} \sff{1+c-b_{2}}{N} \sff{1+c-b_{3}}{N} \sff{1+c-b_{1}-b_{2}-b_{3}}{N}}{\sff{1+c}{N} \sff{1+c-b_{1}-b_{2}}{N} \sff{1+c-b_{1}-b_{3}}{N} \sff{1+c-b_{2}-b_{3}}{N}}}
\beqn{= \ \pfq{9}{8}{c, 1+\z{c}{2}, b_{1}, b_{2}, b_{3}, f-a_{1}-d_{2}-1, f-a_{1}-d_{1}-1, f-d_{1}-d_{2}, -N }{\z{c}{2}, 1+c-b_{1}, 1+c-b_{2}, 1+c-b_{3}, f-d_{1}, f-d_{2}, f-a_{1}-1, 1+c+N}{1}}{2.12}

The (2.12) is equivalent to the Bailey's $_{9}F_8$ transformation ~\cite[ p. 71, (2.4.4.1)]{slt3} with the values of $c$ and $f$ given with (2.7).\\ 

Further, the Bailey's $_{9}F_8(1)$ transformation also follows from the second form of the generalization by selecting $p=a_{1}$. For $p=a_{1}$, it can be verified that  $\lambda =1+a_{1}-f+\z{c}{2}$, $\z{c}{2}+\lambda = 1+a_{1}-f+c=f-d_{1}-d_{2}-1$ and $\z{c}{2}-\lambda = f-a_{1}-1$ and therefore the $_{11}F_{10}$ in the last term of (2.7) converts into $_{9}F_8$ to give (2.12).

\section{Derivations of (2.1), (2.2), (2.4) and (2.7) } 

To derive (2.1), we select 

\beq{e = f + 1,\,\,c = a - f + 1{\mbox{  or  }}f = 1 + a - c}

in the Whipple's transformation (1.4) and observe that the $_{4}F_3(1)$ on the right can be summed by the extended \slz \ theorem (1.5) and after some  simplifications, we get (2.1).

To obtain (2.2), we select
$$ \ez{\alpha}{\sff{f-1,\z{f+1}{2}, a_{1} , a_{2} , f-p ,p+1 }{r} (-1)^{r}}{\sff{\z{f-1}{2}, f-a_{1}, f-a_{2}, p, f-p-1 }{r} r!},\efrc{u}{1}{r!},$$ 
$$ \ez{v}{1}{\sff{f}{r}}, 
 \efrc{\delta}{\sff{d_{1}, d_{2}, -N}{r}}{\sff{1+d_{1}+d_{2}-f-N}{r}},$$\\ 
 in (1.1) and (1.2) and making use of (2.1) and \slz \ theorem, respectively, we get
$$\efnc{\beta}{\sff{f-a_{1}-a_{2}-1, h+1}{n}}{\sff{f-a_{1} , f-a_{2}, h}{n} n!}$$ where the value of $h$ will be as given in (2.3) and $$\efnc{\gamma}{\sff{f-d_{1}, f-d_{2}}{N} }{\sff{f, f-d_{1}-d_{2}}{N}} \z{\sff{d_{1}, d_{2}, -N}{n}(-1)^{n}}{\sff{f-d_{1}, f-d_{2}, f+N}{n}}.$$ \\
Putting these values in (1.3), we obtain the result (2.2).

To find (2.4), we select $a_{2}=2f-2-d_{1}-d_{2}-a_{1}+N$, so that the $_{5}F_4(1)$ on the right of (2.4) converts into a $_{4}F_3(1)$, which can be summed by the extended \slz \ theorem (1.5) and after some simplifications, we get (2.4).

To derive the first form of the generalization of Bailey's $_{9}F_8(1)$ transformation
~\cite[ p. 71, (2.4.4.1)]{slt3}(2.7), we select 
$$ \ez{\alpha}{\sff{f-1,\z{f+1}{2}, a_{1} , d_{1} , d_{2} , f-p ,p+1 }{r} }{\sff{\z{f-1}{2}, f-a_{1}, f-d_{1}, f-d_{2}, p, f-p+1 }{r} r!},\efrc{u}{\sff{f-a_{1}-d_{1}-d_{2}-1}{r}}{r!},$$ 
$$ \ez{v}{\sff{2f-2-d_{1}-d_{2}-a_{1}}{r}}{\sff{f}{r}}, 
 \efrc{\delta}{\sff{1+\z{c}{2}, b_{1}, b_{2}, b_{3}, -N}{r}}{\sff{\z{c}{2}, 1+c-b_{1}, 1+c-b_{2}, 1+c-b_{3}, 1+c+N  }{r}},$$\\ 
 in (1.1) and (1.2), with the restriction that $c=2f-2-d_{1}-d_{2}-a_{1}$ and $3f = 2+b_{1}+b_{2}+b_{3}-d_{1}-d_{2}-a_{1}-N$ and making use of (2.4) and the Dougall's sum, respectively, we get
$$\efnc{\beta}{\sff{f-d_{1}-d_{2}, f-a_{1}-d_{1}-1, f-a_{1}-d_{2}-1, k+1}{n}}{\sff{f-a_{1} , f-d_{2}, f-d_{1}, k}{n} n!}$$ 
and 
$$\efnc{\gamma}{\sff{1+c, 1+c-b_{1}-b_{2}, 1+c-b_{1}-b_{3}, 1+c-b_{2}-b_{3}}{N} }{\sff{1+c-b_{1}, 1+c-b_{2},1+c-b_{3}  1+c-b_{1}-b_{2}-b_{3}}{N}} $$ $$ . \z{\sff{b_{1}, b_{2}, b_{3}, -N}{n}}{\sff{f-b_{1}, f-b_{2}, f-b_{3}, f+N}{n}}.$$ \\

where $k$ needs special attention. Actually, $k=\z{D(A+n)}{B+n}$ with $A=f-d_{1}-d_{2}-1$ and the values of $B$, $D$ and $A'$ are as given in (2.8), (2.9) and (2.10)  will give:

\beq{\z{\sff{k+1}{n}}{\sff{k}{n}}=1+\z{B}{D A} \z{\sff{B+1}{n}}{\sff{B}{n}} \z{\sff{A}{n}}{\sff{A+1}{n}} n}

Putting these values in (1.3) and then using the trick used to prove (2.11), we obtain the first form of the generalization of Bailey's $_{9}F_8(1)$ transformation
~\cite[ p. 71, (2.4.4.1)]{slt3} or the equality of first and middle terms in (2.7).\\ 

To obtain the second form of the generalization of Bailey's $_{9}F_8(1)$ transformation
~\cite[ p. 71, (2.4.4.1)]{slt3}, we proceed as in the proof of the first form, but manipulate the terms containing $k$ as follows:

\beq{\z{\sff{k+1}{n}}{\sff{k}{n}}=\z{\sff{A}{n}}{\sff{A+1}{n}}. \z{\sff{1+\z{c}{2}+\lambda}{n}\sff{1+\z{c}{2}-\lambda}{n}}{\sff{\z{c}{2}-\lambda}{n} \sff{\z{c}{2}+\lambda}{n}}} 
where the value of $\lambda$ is same as given with (2.7). With this value of the terms containing $k$ and the previously calculated values of $\aph_{r}$, $\delta_{r}$, $\beta_{n}$ and $\gamma_{n}$, the (1.3) gives the second form of the generalization of Bailey's $_{9}F_8(1)$ transformation
~\cite[ p. 71, (2.4.4.1)]{slt3}  or the equality of first and last terms in (2.7).

\section{Particular cases of (2.1), (2.2), (2.4) and (2.7) }

Taking $N \to \infty$ in (2.4), we get the extension of a non-terminating $_{5}F_4(1)$ summation \nfff \ as:

\beq{\pfq{7}{6}{f-1,\z{f+1}{2}, a_{1} , f-p ,p+1, d_{1}, d_{2} }{\z{f-1}{2}, f-a_{1} , p ,f-p-1, f-d_{1}, f-d_{2}}{1} \hspace{3.5cm} }
\beqn {= \z{\bgm{f-d_{1}} \bgm{f-d_{2}} \bgm{f-a_{1}} \bgm{f-a_{1}-d_{1}-d_{2}-1}}{\bgm{f} \bgm{f-d_{1}-d_{2}} \bgm{f-a_{1}-d_{1}-1} \bgm{f-a_{1}-d_{2}-1}} \z{\bgm{k}}{\bgm{k+1}}}{4.1}
\beq{\mbox{ provided } Re(f-a_{1}-d_{1}-d_{2}-1)>0,}

where $k$ is given by (2.5) and $h=\z{p(f-p-1)}{a_{1}}$.

Taking $d_{1} \to \infty$ in (4.1), we get the extension of a non-terminating $_{4}F_3(-1)$ summation \nmfft \ as:

\beq{\pfq{6}{5}{f-1,\z{f+1}{2},  f-p ,p+1, a_{1}, d_{2}, }{\z{f-1}{2},  p ,f-p-1, f-a_{1}, f-d_{2}}{-1} \hspace{3.5cm} }
\beqn {= \z{\bgm{f-a_{1}} \bgm{f-d_{2}} }{\bgm{f} \bgm{f-a_{1}-d_{2}-1} } \z{\bgm{k}}{\bgm{k+1}}}{4.2}
\beq{\mbox{ provided } Re(f-2a_{1}-2d_{2}-2)>-1,}

where $k=\z{h(1+d_{2}+a_{1})}{d_{2}-h}$ and $h=\z{p(f-p-1)}{a_{1}}$.

Taking $d_{1}=f/2$ in (4.1), we get the extension of a non-terminating $_{4}F_3(1)$ summation \nfft \ as:

\beq{\pfq{6}{5}{f-1,\z{f+1}{2}, a_{1} , f-p ,p+1, d_{2} }{\z{f-1}{2}, f-a_{1} , p ,f-p-1,  f-d_{2}}{1} \hspace{3.5cm} }
\beqn {= \z{\bgm{f/2} \bgm{f-d_{2}} \bgm{f-a_{1}} \bgm{f/2-a_{1}-d_{2}-1}}{\bgm{f} \bgm{f/2-d_{2}} \bgm{f/2-a_{1}-1} \bgm{f-a_{1}-d_{2}-1}} \z{\bgm{k}}{\bgm{k+1}}}{4.3}
\beq{\mbox{ provided } Re(f-2a_{1}-2d_{2}-2)>0,}

\beqn{\mbox{  where  } k=\z{h(1+a_{1}-f/2)(1+d_{2}+a_{1}-f)}{(d_{2}f/2)-h(1+d_{2}+a_{1}-f/2)} \mbox{ and } h=\z{p(f-p-1)}{a_{1}}}{4.4}

Taking $d_{2} \to \infty$ in (4.3), we get the extension of a non-terminating $_{3}F_2(-1)$ summation ~\cite[ p. 245, (III. 21)]{slt3} \ as:

\beq{\pfq{5}{4}{f-1,\z{f+1}{2}, a_{1} , f-p ,p+1 }{\z{f-1}{2}, f-a_{1} , p ,f-p-1}{-1} \hspace{3.5cm} }
\beqn {= \z{\bgm{f/2}  \bgm{f-a_{1}} }{\bgm{f}  \bgm{f/2-a_{1}-1} } \z{\bgm{k}}{\bgm{k+1}}}{4.5}
\beq{\mbox{ provided } Re(a_{1})>\frac{1}{2},}

\beqn{\mbox{  where  } k=\z{h(1+a_{1}-f/2)}{(f/2)-h} \mbox{ and } h=\z{p(f-p-1)}{a_{1}}}{4.6}

Taking $d_{2}=(f-1)/2$ in (4.1), we get the extension of Dixon's $_{3}F_2(1)$ summation \dixon \ as:

\beq{\pfq{5}{4}{f-1, a_{1} , f-p ,p+1, d_{1} }{ f-a_{1} , p ,f-p-1, f-d_{1}}{1} \hspace{3.5cm} }
\beqn {= \z{\bgm{f-d_{1}} \bgm{\frac{1}{2}f+\frac{1}{2}} \bgm{f-a_{1}} \bgm{\frac{1}{2}f-a_{1}-d_{1}-\frac{1}{2}}}{\bgm{f} \bgm{\frac{1}{2}f-d_{1}+\frac{1}{2}} \bgm{f-a_{1}-d_{1}-1} \bgm{\frac{1}{2}f-a_{1}-\frac{1}{2}}} \z{\bgm{k}}{\bgm{k+1}}}{4.7}
\beq{\mbox{ provided } Re(f-2a_{1}-2d_{1}-1)>0,}

\beqn{\mbox{  where  } k=\z{h(1+d_{1}+a_{1}-f)(\frac{1}{2}+a_{1}-\frac{1}{2}f)}{d_{1}(\frac{1}{2}f-\frac{1}{2}) -h(\frac{1}{2}+d_{1}+a_{1}-\frac{1}{2}f)} \mbox{ and } h=\z{p(f-p-1)}{a_{1}}}{4.8}

Taking $p=\frac{f}{2}$ in (2.4), we get a new very well-poised $_{8}F_7(1)$ summation theorem as:

\beq{\pfq{8}{7}{f-1,\z{f+1}{2}, a_{1} , \frac{f}{2}+1, d_{1}, d_{2}, 2f-2-d_{1}-d_{2}-a_{1}+N, -N}{\z{f-1}{2}, f-a_{1} , \frac{f}{2}-1, f-d_{1}, f-d_{2}, 2+a_{1}+d_{1}+d_{2}-f-N, f+N}{1} \hspace{3.5cm} }
\beqn {= \z{\sff{f}{N} \sff{f-d_{1}-d_{2}}{N} \sff{f-a_{1}-d_{1}-1}{N} \sff{f-a_{1}-d_{2}-1}{N}}{\sff{f-d_{1}}{N} \sff{f-d_{2}}{N} \sff{f-a_{1}}{N} \sff{f-a_{1}-d_{1}-d_{2}-1}{N}} \z{\sff{k+1}{N}}{\sff{k}{N}}}{4.9}

\beqn{\mbox{  where  } k=\z{h(1+d_{1}+a_{1}-f)(1+d_{2}+a_{1}-f)}{d_{1}d_{2}-h(1+d_{1}+d_{2}+a_{1}-f)}} {4.10}

\beqn{\mbox{  and  } h=\z{(\frac{f}{2})(f-1-d_{1}-d_{2}+N)(1-\frac{f}{2})}{p(\frac{f}{2}-1)-a_{1}(2f-2-d_{1}-d_{2}-a_{1}+N)}} {4.11}  

Further, taking $N\to \infty$, we get a new non-terminating $_{6}F_5(1)$ summation theorem for a specific very well-poised hypergeometric function as:

\beq{\pfq{6}{5}{f-1,\z{f+1}{2}, a_{1} , \frac{f}{2}+1, d_{1}, d_{2}} {\z{f-1}{2}, f-a_{1} , \frac{f}{2}-1, f-d_{1}, f-d_{2}}{1} \hspace{3.5cm} }
\beqn {= \z{\bgm{f-d_{1}} \bgm{f-d_{2}} \bgm{f-a_{1}} \bgm{f-a_{1}-d_{1}-d_{2}-1}}{\bgm{f} \bgm{f-d_{1}-d_{2}} \bgm{f-a_{1}-d_{1}-1} \bgm{f-a_{1}-d_{2}-1}} \z{\bgm{k}}{\bgm{k+1}}}{4.12}
\beq{\mbox{ provided } Re(f-a_{1}-d_{1}-d_{2}-1)>0,}

\beqn{\mbox{  where  } k=\z{h(1+d_{1}+a_{1}-f)(1+d_{2}+a_{1}-f)}{d_{1}d_{2}-h(1+d_{1}+d_{2}+a_{1}-f)}} {4.13}

\beqn{\mbox{  and  } h=\z{(\frac{f}{2})(1-\frac{f}{2})}{-a_{1}}} {4.14}  

Taking $p=\frac{f}{2}+\frac{1}{2}$ in (2.4), we get a new summation theorem contiguous to Dougall's $_{7}F_6(1)$ theorem from very well-poised pair as:

\beq{\pfq{7}{6}{f-1,\frac{f}{2}+\frac{3}{2}, a_{1} , d_{1}, d_{2}, 2f-2-d_{1}-d_{2}-a_{1}+N, -N}{\frac{f}{2}-\frac{3}{2}, f-a_{1} , f-d_{1}, f-d_{2}, 2+a_{1}+d_{1}+d_{2}-f-N, f+N}{1} \hspace{3.5cm} }
\beqn {= \z{\sff{f}{N} \sff{f-d_{1}-d_{2}}{N} \sff{f-a_{1}-d_{1}-1}{N} \sff{f-a_{1}-d_{2}-1}{N}}{\sff{f-d_{1}}{N} \sff{f-d_{2}}{N} \sff{f-a_{1}}{N} \sff{f-a_{1}-d_{1}-d_{2}-1}{N}} \z{\sff{k+1}{N}}{\sff{k}{N}}}{4.15}

\beqn{\mbox{  where  } k=\z{h(1+d_{1}+a_{1}-f)(1+d_{2}+a_{1}-f)}{d_{1}d_{2}-h(1+d_{1}+d_{2}+a_{1}-f)}} {4.16}

\beqn{\mbox{  and  } h=\z{(\frac{f}{2}+\frac{1}{2})(f-1-d_{1}-d_{2}+N)(\frac{3}{2}-\frac{f}{2})} {(\frac{f}{2}+\frac{1}{2}) (\frac{f}{2}-\frac{3}{2})-a_{1}(2f-2-d_{1}-d_{2}-a_{1}+N)}} {4.17}  

Taking $N\to \infty$ in (4.15), we get a new summation theorem contiguous to classical $_{5}F_4(1)$ from very well-poised pair as:

\beq{\pfq{5}{4}{f-1,\frac{f}{2}+\frac{3}{2}, a_{1} , d_{1}, d_{2}} {\frac{f}{2}-\frac{3}{2}, f-a_{1} , f-d_{1}, f-d_{2}}{1} \hspace{3.5cm} }
\beqn {= \z{\bgm{f-d_{1}} \bgm{f-d_{2}} \bgm{f-a_{1}} \bgm{f-a_{1}-d_{1}-d_{2}-1}}{\bgm{f} \bgm{f-d_{1}-d_{2}} \bgm{f-a_{1}-d_{1}-1} \bgm{f-a_{1}-d_{2}-1}} \z{\bgm{k}}{\bgm{k+1}}}{4.18}
\beq{\mbox{ provided } Re(f-a_{1}-d_{1}-d_{2}-1)>0,}

\beqn{\mbox{  where  } k=\z{h(1+d_{1}+a_{1}-f)(1+d_{2}+a_{1}-f)}{d_{1}d_{2}-h(1+d_{1}+d_{2}+a_{1}-f)}} {4.19}

\beqn{\mbox{  and  } h=\z{(\frac{f}{2}+\frac{1}{2})(\frac{f}{2}-\frac{3}{2})}{a_{1}}} {4.20}  

Taking $p\to \infty$ in (2.4), we get a new summation theorem contiguous to Dougall's $_{7}F_6(1)$ theorem from second last numerator parameter as:

\beq{\pfq{7}{6}{f-1,\z{f+1}{2}, a_{1} , d_{1}, d_{2}, 2f-2-d_{1}-d_{2}-a_{1}+N, -N}{\z{f-1}{2}, f-a_{1} , f-d_{1}, f-d_{2}, 2+a_{1}+d_{1}+d_{2}-f-N, f+N}{1} \hspace{3.5cm} }
\beqn {= \z{\sff{f}{N} \sff{f-d_{1}-d_{2}}{N} \sff{f-a_{1}-d_{1}-1}{N} \sff{f-a_{1}-d_{2}-1}{N}}{\sff{f-d_{1}}{N} \sff{f-d_{2}}{N} \sff{f-a_{1}}{N} \sff{f-a_{1}-d_{1}-d_{2}-1}{N}} \z{\sff{k+1}{N}}{\sff{k}{N}}}{4.21}

\beqn{\mbox{  where  } k=\z{h(1+d_{1}+a_{1}-f)(1+d_{2}+a_{1}-f)}{d_{1}d_{2}-h(1+d_{1}+d_{2}+a_{1}-f)}} {4.22}

\beqn{\mbox{  and  } h=(1-f+d_{1}+d_{2}-N)} {4.23}  

which also leads to classical $_{5}F_4(1)$ summation theorem \nfff \ as $N\to \infty$. \\ 
Further, on the above mentioned lines, the particular cases of (2.2) and (2.7) can also be obtained, which will produce new hypergeometric transformations.

\section{What next?}
This work has filled gap in the existing literature by adding a number of new hypergeometric summations and transformations in the form of new results and extensions of classical hypergeometric theorems to the Slater's list ~\cite[ p. 243-245, Appendix III)]{slt3} and will, certainly, boost the investigations of many more general and significant hypergeometric summations and transformations. Additional results are under preparation and will be announced soon.

\begin{flushleft}
Department of Mathematics,\\ 
School of Enigineering,\\
Sir Padampat Singhania University, \\
Bhatewar, Udaipur-313 601, Rajasthan, INDIA.\\

E-Mail : {\it yashoverdhan.vyas@spsu.ac.in} \\
\end{flushleft}
\end{document}